\documentclass[fleqn]{mat01}
\usepackage{times,mathtimy,amssymb,latexsym,boxit}
\begin{document}

\setcounter{page}{213}
\firstpage{213}

\newtheorem{theo}{Theorem}
\renewcommand\thetheo{\arabic{theo}}
\newtheorem{theor}[theo]{\bf Theorem}
\newtheorem{lem}{Lemma}
\newtheorem{coro}{\rm COROLLARY}
\def\notat{\trivlist\item[\hskip\labelsep{\it Notation.}]}
\def\noot{\trivlist\item[\hskip\labelsep{\it Note.}]}
\def\pot{\trivlist\item[\hskip\labelsep{\it Proofs of Theorems $1$ and $2$.}]}

\title{Analogues of Euler and Poisson summation formulae}

\markboth{Vivek~V Rane}{Analogues of Euler and Poisson summation formulae}

\author{VIVEK~V RANE}

\address{Department of Mathematics, The Institute of Science, 15, Madam
Cama Road, Mumbai~400~032, India}

\volume{113}

\mon{August}

\parts{3}

\Date{MS received 2 February 2002; revised 20 January 2003}

\begin{abstract}
Euler--Maclaurin and Poisson analogues of the summations $\sum_{a < n
\leq b} \chi(n) f(n)$, $\sum_{a < n \leq b} d(n) f(n)$, $\sum_{a < n \leq b}
d(n) \chi (n) f(n)$ have been obtained in a unified manner, where $(\chi
(n))$ is a periodic complex sequence; $d(n)$ is the divisor function and
$f(x)$ is a sufficiently smooth function on $[a,b]$. We also state a
generalised Abel's summation formula, generalised Euler's summation
formula and Euler's summation formula in several variables.
\end{abstract}


\maketitle

\section{Introduction}

Voronoi \cite{4} conjectured that if $(c(n))$ is a given arithmetical
function and if $f$ is a continuous function on an interval $[a,b]$ with
only finite number of maxima and minima there, then there exist analytic
functions $\alpha(x)$ and $\delta (x)$, depending only upon $c(n)$ (and
not upon $f(x)$) such that
\begin{equation*}
\sum\limits_{a \leq n \leq b}\!\!\!{'}\ c(n) f(n) = \int_{a}^{b} f(x) \delta (x)
{\rm d}x + \sum\limits_{n = 1}^{\infty} c(n) \int_{a}^{b} f(x) \alpha
(nx) {\rm d}x.
\end{equation*}
The prime on the summation sign $\sum{'}_{a \leq n \leq b}\ c(n) f(n)$
means that if $n = a$ or $n = b$, only $(1/2) c(a) f(a)$ or $(1/2) c(b)
f(b)$ respectively is counted. In the special case $c(n) = d(n)$, where
$d(n) = $ the number of divisors of $n$, Voronoi obtained the formula
\begin{equation*}
\sum\limits_{a \leq n \leq b}\!\!\!{'}\ d(n) f(n) = \int_{a}^{b} (\log\ x + 2\gamma)
f(x) {\rm d}x + \sum\limits_{n = 1}^{\infty} d(n) \int_{a}^{b} f(x)
\alpha (nx) {\rm d}x,
\end{equation*}
with $\alpha (x) = 4 K_{0} (4\pi \sqrt{x}) - 2\pi Y_{0} (4\pi
\sqrt{x})$, where $K_{0}$ and $Y_{0}$ are the well-known Bessel
functions and $\gamma$ is the Euler's constant. Obviously, Voronoi's
summation formula is a generalisation of Poisson's summation formula.
Berndt \cite{1} and Berndt and Schoenfeld \cite{2} have given the Euler--Maclaurin
and the Poisson analogues of the summation $\sum_{a < n \leq b} \chi (n)
f(n)$, where $(\chi(n))$ is a periodic sequence of complex numbers.

The object of this paper is as follows:
\begin{enumerate}
\renewcommand{\labelenumi}{(\arabic{enumi})}
\item To give the analogues of the Euler--Maclaurin's summation formula
of $\sum_{a\!<n \leq b}\! \chi(n)\! f(n)$, $\sum_{a < n \leq b} d(n) f(n)$ and
$\sum_{a < n \leq b} d(n) \chi(n) f(n)$, where $f(x)$ is a sufficiently
smooth function on the interval $[a, b]; d(n)$ is the divisor function
and $(\chi (n))$ is a periodic sequence of complex numbers of period
$k$ (see Theorem~1).

\item To give the analogues of the Poisson's summation formula for
$\sum_{a < n \leq b} \chi(n) f(n)$, $\sum_{a < n \leq b} d(n) f(n)$ and
$\sum_{a < n \leq b} d(n) \chi(n) f(n)$ (see Theorem~2). These
will be accomplished virtually effortlessly without the use of complex
contour integration. Incidentally in particular, choosing $\chi(n)$ to
be the constant sequence $(1)_{n}$ so that $k = 1$, we can obtain
Euler--Maclaurin and Poisson summation formulae for $\sum_{a < n \leq b}
f(n)$ from the corresponding results for $\sum_{a < n \leq b} \chi (n)
f(n)$.
\end{enumerate}

Next, we introduce our notation and state our theorems.

\begin{notat}
In what follows, the summation $\sum_{n =-\infty}^{\infty}$ or
$\sum_{n}$ will always mean the summation in the sense $\lim_{N
\rightarrow \infty} \sum_{|n| \leq N}$; and $\sum{'}_{n=-\infty}^{\infty}$
or $\sum{'}_{n}$ shall denote exclusion of the term corresponding to $n
=0$. For $x$ real, $[x]$ denotes the integral part of $x$. We shall
write $\psi(x) = x - [x] -(1/2)$. Thus for non-integer $x$, we have
$\psi(x) =-\sum{'}_{n} ({\rm e}^{2\pi in x}/2\pi in)$. We shall call a
Riemann-integrable function $f$ on an interval $[a,b]$ of the real line,
a good function on $[a,b]$, if it admits the interchange of $\sum$ and
$\int$ in the Riemann--Stieltjes integral $\int_{a}^{b} \psi (\alpha x +
\beta) {\rm d}f(x)$, i.e.,
\begin{equation*}
-\int_{a}^{b} \sum\limits_{n=-\infty}^{\infty}\!\!\!{'}\
\frac{{\rm e}^{2\pi in (\alpha x + \beta)}}{2\pi in} {\rm
d}f (x) =-\sum_{n=-\infty}^{\infty}\!\!\!{'}\ \int_{a}^{b} \frac{{\rm e}^{2\pi
in(\alpha x +\beta)}}{2\pi in} {\rm d}f (x),
\end{equation*}
where $\alpha, \beta$ are some real constants with $0 < \alpha \leq 1$.
 (If $f$ is differentiable with its derivative finite
and integrable, whether in the Riemann sense or in the Lebesgue sense,
on $[a,b]$, then $f$ is good on $[a,b]$. Also, if $f$ is a function of
bounded variation, then $f$ is good.)

For an integer $r \geq 0$, we shall write
\begin{equation*}
\psi_{r} (x) =-\sum_{n}{'}\ \frac{{\rm e}^{2\pi in x}}{(2\pi
in)^{r +1}} = \frac{B_{r+1} (x)}{(r + 1)!},
\end{equation*}
where $B_{r} (x)$ is the Bernoulli polynomial of degree $r$ and $(r +
1)!$ denotes $(r + 1)$ factorial. Note that $\psi_{0} (x) = \psi (x)$
and $({\rm d}/{\rm d}x) \psi_{r} (x) = \psi_{r-1} (x)$ for $r \geq 1$.
\end{notat}

Next, we state our theorems.

\begin{theor}[\!]
Let the function $f$ have continuous derivatives up to $(R + 1)$th order
on $[a,b]$.

\begin{enumerate}
\renewcommand{\labelenumi}{\rm (\Roman{enumi})}
\leftskip .5pc
\item Let $(\chi(n))$ be a complex sequence with period $k$ and let
$d(n)$ be the divisor function. Then{\rm ,} we have
\begin{align*}
\hskip -1.3pc \sum\limits_{a < n \leq b} \chi(n) f(n) &= \frac{1}{k}
\sum\limits_{l=1}^{k} \chi(l) \int_{a}^{b} f(u){\rm d}u\\
&\quad\ - \sum\limits_{r=0}^{R} (-k)^{r} \left\lbrace f^{(r)} (b)
\sum\limits_{n}{'}\ \tau(\chi, -n) \frac{{\rm e}^{2\pi in
b/k}}{(2\pi in)^{r+1}} \right.\\
&\quad\ \left. - f^{(r)} (a) \sum\limits_{n}{'}\ \tau
(\chi, -n) \frac{{\rm e}^{2\pi in a/k}}{(2\pi in)^{r+1}} \right\rbrace\\
&\quad\ + (-k)^{R} \sum\limits_{n}{'}\ \! \frac{\tau (\chi, -n)}{(2\pi
in)^{R+1}} \int_{a}^{b}\!\! f^{(R + 1)}(u) {\rm e}^{2\pi in u/k} {\rm
d}u,\hskip 1.5pc
\end{align*}
where $\tau (\chi, n) = \sum_{l=1}^{k} \chi (l) {\rm e}^{2\pi in l/k}$.

\item We have\vspace{-.4pc}
\begin{align*}
\hskip -3pc &\sum\limits_{a < n \leq b} d(n) f(n)\\
\hskip -3pc &\quad\ = \int_{a}^{b} f(u) {\rm
d}(u) \left( \sum\limits_{m \leq b} \frac{1}{m}\right)\\
\hskip -3pc &\qquad\ + \sum\limits_{r=0}^{R} (-1)^{(r+1)} \left\lbrace f^{(r)} (b)
\sum\limits_{n}{'}\frac{1}{(2\pi in)^{r + 1}} \left( \sum\limits_{m \leq
b} m^{r} {\rm e}^{2\pi in b/m} \right) \right.\\
\hskip -3pc &\qquad \ \left. - f^{(r)} (a) \sum\limits_{n}{'} \frac{1}{(2\pi in)^{r+1}}
\left( \sum\limits_{m \leq b} m^{r} {\rm e}^{2\pi in a/m} \right)
\right\rbrace\\
\hskip -3pc &\qquad\ + (-1)^{R} \sum\limits_{n}{'} \frac{1}{(2\pi in)^{R+1}}
\int_{a}^{b} f^{(R+1)} (u) \left( \sum\limits_{m \leq b} m^{R} {\rm
e}^{2\pi in u/m} \right) {\rm d}u.
\end{align*}

\item We have\vspace{-.4pc}
\begin{align*}
\hskip -3pc \sum\limits_{a < n \leq b} \chi(n) d(n) f(n) &= \frac{1}{k}
\sum\limits_{r_{1} = 1}^{k} \sum\limits_{r_{2} = 1}^{k} \chi (r_{1} r_{2})
\left( \sum\limits_{\substack{m \leq b\\ m \equiv r_{1} ({\rm mod}\
k)}} \frac{1}{m} \int_{a}^{b} f(u) {\rm d}u \right)\\
&\quad\ + \sum\limits_{r = 0}^{R} (-1)^{r+1} k^{r} \sum\limits_{r_{1} =
1}^{k} \sum\limits_{r_{2} = 1}^{k} \chi (r_{1} r_{2})\\
\hskip -3pc &\quad\ \times \left\lbrace
f^{(r)} (b) \sum_{n}{'} \frac{{\rm e}^{-2\pi in r_{2}/k}}{(2\pi
in)^{r+1}} \cdot \sum\limits_{\substack{m \leq b\\ m \equiv r_{1} ({\rm
mod}\ k)}} {\rm e}^{2\pi in b/mk} \right.\\
\hskip -3pc &\qquad\quad\ \left. -f^{(r)} (a) \sum\limits_{n}{'} \frac{{\rm e}^{-2\pi in
r_{2}/k}}{(2\pi in)^{r + 1}} \cdot \sum\limits_{\substack{m \leq b\\ m
\equiv r_{1} ({\rm mod}\ k)}} {\rm e}^{2\pi in a/mk} \right\rbrace\\
\hskip -3pc &\quad\ + (-k)^{R} \sum\limits_{r_{1} = 1}^{k} \sum\limits_{r_{2} =
1}^{k} \chi (r_{1} r_{2}) \sum\limits_{n}{'} {\rm e}^{-2\pi in
r_{2}/k} \int_{a}^{b} f^{(R + 1)} (u)\\
\hskip -3pc &\qquad\ \left( \sum\limits_{\substack{m
\leq b\\ m \equiv r_{1}({\rm mod}\ k)}} {\rm e}^{2\pi in u/mk} \right)
{\rm d}u.
\end{align*}
\end{enumerate}
\end{theor}

\begin{theor}[\!]
Let a function $f$ be good on the interval $[a,b]$. Let $(\chi (n))$ be
a complex sequence of period $k$ and let $d(n)$ be the divisor function.
Then{\rm ,} we have

\begin{enumerate}
\renewcommand{\labelenumi}{\rm (\Roman{enumi})}

\leftskip .5pc
\item{}$\left.\right.$\vspace{-3pc}

\begin{equation*}
\hskip -2pc \sum\limits_{a < n \leq b} \chi(n) f(n) = \frac{1}{k} \sum\limits_{n =-
\infty}^{\infty} \tau (\chi, n) \cdot \int_{a}^{b} f(u) {\rm e}^{-2\pi
in u/k} {\rm d}u,
\end{equation*}
where $\tau (\chi, n) = \sum_{l =1}^{k} \chi(l) {\rm e}^{2\pi
in l/k}$.

\item We have\vspace{-.5pc}
\begin{equation*}
\hskip -1.5pc\sum\limits_{a < n \leq b} d(n) f(n) = \sum\limits_{n=-\infty}^{\infty}
\int_{a}^{b} f(u) \left( \sum\limits_{m \leq b} \frac{1}{m} {\rm
e}^{2\pi in u/m} \right) {\rm d}u.
\end{equation*}
\item \vspace{-.5pc}
\begin{align*}
\hskip -1.7pc &\sum\limits_{a < n \leq b} \chi(n) d(n) f(n)\\
\hskip -1.7pc &\quad\ = \sum\limits_{r_{1} =
1}^{k} \sum\limits_{r_{2} = 1}^{k} \chi(r_{1} r_{2}) \sum\limits_{n=-
\infty}^{\infty}\!\!\!{'}\ \ {\rm e}^{-2\pi in r_{2}/k}\\
\hskip -1.7pc &\qquad\ \times \int_{a}^{b} {\rm d}uf (u)
 \left( \sum\limits_{\substack{m \leq b\\ m \equiv r_{1} ({\rm mod}\ k)}}
\frac{1}{m} {\rm e}^{2\pi in u/mk} \right)\\
\hskip -1.7pc &\qquad\ + \frac{1}{k}
\sum\limits_{r_{1} = 1}^{k} \sum\limits_{r_{2} = 1}^{k} \chi (r_{1}
r_{2}) \left( \int_{a}^{b} f(u) {\rm d}u \right) \left(
\sum\limits_{\substack{m \leq b\\ m \equiv r_{1} ({\rm mod}\ k)}} \frac{1}{m}
\right).
\end{align*}
\end{enumerate}\vspace{-.7pc}
\end{theor}

Berndt \cite{1} and Berndt and Schoenfeld \cite{2} have given the Euler--Maclaurin
and the Poisson analogues for the summation $\sum_{a < n \leq b} \chi (n)
f(n)$ by different methods. Their results are similar to our results for
$\sum_{a < n \leq b} \chi (n) f(n)$. Jutila \cite{3} obtained
transformation formulae in analytic number theory, where he deals with
the summations $\sum_{n} d(n)f(n)$ and $\sum_{n} a(n) f(n); a(n)$ being
the $n$th Fourier coefficient of a cusp form. We can also obtain the
Euler--Maclaurin and the Poisson analogues of $\sum_{n} a(n) f(n)$,
using our approach. More generally, we can deal with the summation
$\sum_{n} r(n)f(n)$, where $r(n)$ is the Fourier coefficient of a
periodic integrable function $g(x)$ of period 1. Thus writing $r(x) =
\int_{0}^{1} g(u) {\rm e}^{-2\pi ix u} {\rm d}u$, we have $r(n) =
\int_{0}^{1} g(u) {\rm e}^{-2\pi in u} {\rm d}u$.

Note that\vspace{-.3pc}
\begin{align*}
{\frac{\rm d}{{{\rm d}}\hbox{$x$}}} r(x) &= (-2\pi ix) \int_{0}^{1} g(u) {\rm e}^{-
2\pi ixu} {\rm d}u,\\
\frac{\rm d^{2}}{{{\rm d}}\hbox{$x$}^{2}} r(x) &= (-2\pi ix)^{2} \int_{0}^{1}
g(u) {\rm e}^{-2\pi ixu} {\rm d}u,
\end{align*}
and so on.

Thus $r(x)$ is a smooth function. In particular, we may choose $g(u) =
h_{y} (u) = h(u + iy)$ for a fixed $y > 0$, where $h(z) = h(x + iy)$ is a
cusp form under consideration as in the case of Jutila \cite{3}. If $f$
is a smooth function on $[a,b]$, then $\phi (x) = f(x) \cdot r(x)$ is a
smooth function on $[a,b]$ and we can apply Euler--Maclaurin or Poisson
summation formula to the summation $\sum_{a < n \leq b} \phi(n)$.

We actually show that all the above results can be obtained using the
two facts, namely,

\begin{enumerate}
\renewcommand{\labelenumi}{(\arabic{enumi})}
\item generalised Euler's summation formula (see Corollary~1 of Lemma~1),

\item $\psi (x) = x-[x] - \frac{1}{2} =-\sum{'}_{n=-\infty}^{\infty} {\rm
e}^{(2\pi inx/2\pi in)}$, the series being convergent boundedly.
\end{enumerate}\vspace{-.5pc}

Next we prove the theorems. For this, we state our main results as
lemmas and derive the proofs of our theorems from these lemmas.

\begin{lem}{\rm (}Generalised Abel's summation formula{\rm )}.\ \ Let $(\lambda
(n))_{n=-\infty}^{n=\infty}$ be a strictly increasing sequence of real
numbers such that $\lambda (n) \rightarrow \infty$ as $n \rightarrow
\infty${\rm ;} and $\lambda(n) \rightarrow -\infty$ as $n \rightarrow -
\infty$. Let $f$ be a function on an interval $[a,b]$ such that $f$ is
continuous from left at every point of the sequence $(\lambda (n))$ with
$a < \lambda(n) \leq b$. Let $(c(n))$ be a complex sequence and let $S(t)
= \sum_{\lambda_{0} < \lambda (n) \leq t} c(n)${\rm ,} where $\lambda_{0}$ is
a fixed constant. Then
\begin{align*}
\sum\limits_{a < \lambda (n) \leq b} c(n) f( \lambda (n)) &= \int_{a}^{b}
f(t) {\rm d}S(t)\\
&= f(b) S(b) - f(a) S(a) - \int_{a}^{b} S(t) {\rm d}f (t).
\end{align*}
\end{lem}

A corresponding result may be given, if $f$ is continuous from right at
every point $\lambda(n)$ with $a < \lambda(n) \leq b$.

\begin{coro}{\rm (}Generalised Euler's summation formula{\rm )}$\left.\right.$\vspace{.5pc}

\noindent Let $f$ be Riemann-integrable on the interval $[a,b]$ such
that $f$ is continuous from left at every integer $n$ with $a < n \leq
b$. Then
\begin{align*}
\sum\limits_{a < n \leq b} f(n) &= \int_{a}^{b} f(u) {\rm d}u +
\int_{a}^{b} \left(u-[u] - \frac{1}{2} \right) {\rm d}f (u)\\
&\quad\ + f(a) \left(a -[a] - \frac{1}{2} \right) -f(b) \left( b-[b] -
\frac{1}{2} \right).
\end{align*}
\end{coro}

\begin{noot}
If $f$ is such that its derivative $f'$ exists and is finite on $[a,b]$
and is integrable on $[a,b]$ (either in Riemann or Lebesgue sense), then
we can replace $\int_{a}^{b} \big(u-[u] - \frac{1}{2} \big) {\rm d}f
(u)$ by $\int_{a}^{b} \big( u-[u] - \frac{1}{2} \big) f' (u) {\rm
d}u$.
\end{noot}\vspace{.5pc}

\begin{coro}{\rm (}Euler's summation formula for two
variables{\rm )}$\left.\right.$\vspace{.5pc}

\noindent Let $f(x,y)$ be a function of two variables such that its
partial derivatives up to second order are continuous in the rectangle
$(a \leq x \leq b, c \leq y \leq d)${\rm ,} where $a, b, c, d$ are integers.
Then with obvious notations{\rm ,}
\begin{align*}
\sum\limits_{c < n \leq d} \sum\limits_{a < m \leq b} f(m, n) &=
\int_{a}^{b} \int_{c}^{d} f(x, y) {\rm d}x {\rm d}y\\
&\quad\ + \int_{a}^{b} \int_{c}^{d} f_{x} (x, y) (x - [x]) {\rm d}x {\rm d}y\\
&\quad\ + \int_{a}^{b} \int_{c}^{d} f_{y} (x, y) (y - [y]) {\rm d}x {\rm d}y\\
&\quad\ + \int_{a}^{b} \int_{c}^{d} f_{xy} (x, y) (x - [x]) (y - [y])
{\rm d}x {\rm d}y.
\end{align*}
\end{coro}

\begin{proof}
Using Euler's summation formula (for fixed $n$), we have
\begin{equation*}
\sum\limits_{a < m \leq b} f(m, n) = \int_{a}^{b} f(x, n) {\rm d}x +
\int_{a}^{b} f_{x} (x, n) (x -[x]) {\rm d}x.
\end{equation*}
Hence,\vspace{-.5pc}
\begin{align*}
\sum\limits_{c < n \leq d} \left( \sum\limits_{a < m \leq b} f(m, n)
\right) &= \int_{a}^{b} {\rm d}x \left( \sum\limits_{c < n \leq d} f(x,
n) \right)\\
&\quad\ + \int_{a}^{b} {\rm d}x (x -[x]) \left( \sum\limits_{c < n \leq
d} f_{x} (x, n) \right).
\end{align*}
Using Euler's summation formula once more, for the summations $\sum_{c <
n \leq d} f(x, n)$ and $\sum_{c < n \leq d} f_{x} (x, n)$, we get the
result as stated.\vspace{.4pc}
\end{proof}

\begin{coro}$\left.\right.$\vspace{.5pc}

\noindent We have for integers $r, k$ with $0 \leq r < k${\rm ,} if $f$ is
continuous from left at the integer $n \equiv r({\rm mod}\ k)$ with $a
< n \leq b${\rm ,} then
\begin{enumerate}
\renewcommand{\labelenumi}{\rm (\Roman{enumi})}
\leftskip .5pc
\item $\left.\right.$\vspace{-2.5pc}

\begin{align*}
\hskip -1pc \sum\limits_{\substack{a < n \leq b\\ n \equiv r({\rm mod}\ k)}}
f(n) &= \frac{1}{k} \int_{a}^{b} f(u) {\rm d}u + \int_{a}^{b} \psi \left(
\frac{u-r}{k} \right) {\rm d}f (u)\\
\hskip -1pc &\quad\ + f(a) \cdot \psi \left( \frac{a -
r}{k} \right) - f(b) \cdot \psi \left( \frac{b - r}{k} \right),
\end{align*}
provided $f$ is Riemann-integrable on $[a, b]$.

\item Putting $r =0$ we get for $m \geq 1$, if $f$ is continuous on $[a,
b]$
\begin{align*}
\hskip -1pc \sum\limits_{a/m < n \leq b/m} f(mn) &= \frac{1}{m} \int_{a}^{b} f(u)
{\rm d}u + \int_{a}^{b} \psi \left( \frac{u}{m} \right) {\rm d}f (u) +
f(a) \psi \left( \frac{a}{m} \right)\\
\hskip -1pc &\quad\ - f(b) \psi \left( \frac{b}{m} \right).
\end{align*}

\item We have for integer $m \geq 1$ and integer $k \geq 1${\rm ,}
\begin{align*}
\hskip -1pc \sum\limits_{\substack{a/m < n \leq b/m\\ n\equiv r({\rm mod}\ k)}}
f(mn) &= \frac{1}{km} \int_{a}^{b} f(u) {\rm d}u + \int_{a}^{b} \psi
\left( \frac{(u/m) - r}{k} \right) {\rm d}f (u)\\
\hskip -1pc &\quad\ + f(a) \psi \left( \frac{(a/m) - r}{k} \right) - f(b) \psi \left(
\frac{(b/m) -r}{k} \right),
\end{align*}
provided $f$ is continuous on $[a,b]$.
\end{enumerate}
\end{coro}

\begin{pot}
Firstly, we shall find the Poisson and the Euler--Maclaurin analogues of
$\sum_{a < n \leq b} \chi(n) f(n)$. We have
\begin{align*}
\sum\limits_{a < n \leq b} \chi(n) f(n) &= \sum\limits_{l = 1}^{k}
\chi(l) \left( \sum\limits_{n \equiv l ({\rm mod}\ k)} f(n) \right)\\
&= \sum\limits_{l = 1}^{k} \chi(l) \left\lbrace \frac{1}{k} \int_{a}^{b}
f(u) {\rm d}u  + \int_{a}^{b} \psi \left( \frac{u - l}{k} \right) {\rm d}f (u)\right.\\
&\quad\ \left. + f(a) \psi \left( \frac{a - l}{k} \right) - f(b) \psi \left( \frac{b -
l}{k} \right) \right\rbrace.
\end{align*}
First, we shall obtain Poisson analogue. Integrating by parts, we have
\begin{align*}
\hskip -2pc \int_{a}^{b} \psi \left( \frac{u - l}{k} \right) {\rm d}f (u) &= -
\int_{a}^{b} \left( \sum\limits_{n =-\infty}^{\infty}\!\!\!{'}\ \ \frac{{\rm
e}^{2\pi in((u- l)/k)}}{2\pi in} \right) {\rm d}f (u)\\
\hskip -2pc &= -\sum\limits_{n=-\infty}^{\infty}\!\!\!{'}\ \ \frac{1}{2\pi in} \int_{a}^{b}
{\rm e}^{2\pi in ((u-l)/k)} {\rm d}f (u)\\
\hskip -2pc &= -\sum\limits_{n}{'} \frac{1}{2\pi in} \bigg\lbrace [{\rm e}^{2\pi
in ((u-l)/k)} f (u)]_{u=a}^{u=b} \\
\hskip -2pc &\quad\ \left. - \frac{2\pi in}{k} \int_{a}^{b} f(u) {\rm e}^{2\pi in((u-
l)/k)} {\rm d}u \right\rbrace\\
\hskip -2pc &= -\sum\limits_{n}{'} \frac{1}{2\pi in} \left( {\rm e}^{2\pi in ((b-
l)/k)} f(b) - {\rm e}^{2\pi in ((a-l)/k)} f(a) \right)\\
\hskip -2pc &\quad\ + \frac{1}{k} \sum\limits_{n}{'} \int_{a}^{b} f(u) {\rm e}^{2\pi in ((u-l)/k)} {\rm d}u.
\end{align*}
Thus Theorem~2(I) follows.

Next, we prove the Euler--Maclaurin analogue of $\sum_{a < n \leq
b} \chi(n) f(n)$. Integrating by parts, we have
\begin{align*}
\hskip -2pc \int_{a}^{b} \psi \left( \frac{u-l}{k} \right) {\rm d}f (u) &=
\int_{a}^{b} \psi \left( \frac{u-l}{k} \right) f' (u) {\rm d}u\\
\hskip -2pc &= k \int_{a}^{b} ({\rm d}/{\rm d}u) \psi_{1} \left( \frac{u -l}{k}
\right) f' (u) {\rm d}u\\
\hskip -2pc &= k \left\lbrack \psi_{1} \left( \frac{u-l}{k} \right) f' (u)
\right\rbrack_{u=a}^{b} -k \int_{a}^{b} \psi_{1} \left( \frac{u-l}{k}
\right) f'' (u) {\rm d}u\\
\hskip -2pc &= k \left( \psi_{1} \left( \frac{b-l}{k} \right) f'(b) - \psi_{1} \left(
\frac{u-l}{k} \right) f'(a) \right)\\
\hskip -2pc &\quad\ - k \int_{a}^{b} \psi_{1} \left(
\frac{u-l}{k} \right) f'' (u) {\rm d}u\\
\hskip -2pc &= k \left( \psi_{1} \left( \frac{b-l}{k} \right) f'(b) - \psi_{1}
\left( \frac{a-l}{k} \right) f' (a) \right)\\
\hskip -2pc &\quad\ - k^{2}\int_{a}^{b} \frac{{\rm d}}{{\rm d}u} \psi_{2} \left(
\frac{u-l}{k} \right) f'' (u) {\rm d}u,
\end{align*}
and so on. Thus, we get
\begin{align*}
\hskip -2pc \sum\limits_{a < n \leq b} \chi(n) f(n) &= \sum\limits_{l=1}^{k} \chi(l)
\left\lbrace \frac{1}{k} \int_{a}^{b} f(u) {\rm d}u\right.\\
&\quad\ + \sum\limits_{r=0}^{R} (-1)^{r+1} k^{r} \left( \psi_{r} \left(
\frac{b\!-\!l}{k} \right) f^{(r)} (b) \!-\! \psi_{r} \left( \frac{a\!-\!l}{k}
\right) f^{(r)} (a) \right)\\
&\quad\ \left. - k^{R} \int_{a}^{b} \psi_{R} \left( \frac{u-l}{k} \right) f^{(R
+ 1)} (u) {\rm d}u \right\rbrace.
\end{align*}
This gives Theorem~1(I). Next, we deal with $\sum_{a < n \leq b} {\rm
d}(n) f(n)$.\vspace{.4pc}

Now $\sum_{a < r \leq b} d(r) f(r) = \sum_{a < mn \leq b} f(mn)$.\vspace{.4pc}
\end{pot}

Noting that the summation $a < mn \leq b$ means summation over lattice
points between rectangular hyperbolae $xy=b$ and $xy=a$, the upper
hyperbola included and the lower hyperbola excluded, we get
\begin{align*}
\sum\limits_{a < r \leq b} d(r) f(r) &= \sum\limits_{m \leq b} \left(
\sum\limits_{a/m < n \leq b/m} f(mn)\right)\\
&= \sum\limits_{m \leq b} \frac{1}{m} \left( \int_{a}^{b} f(u) {\rm d}u
\right) + \sum\limits_{m \leq b} \left( \int_{a}^{b} \psi \left(
\frac{u}{m} \right) {\rm d}f (u) \right.\\
&\quad\ \left. + f(a) \psi \left( \frac{a}{m} \right)
- f(b) \psi \left( \frac{b}{m} \right) \right).
\end{align*}
Substituting the series for $\psi (u/m)$, we can obtain both the Poisson
and the Euler--Maclaurin analogues for $\sum_{a < n \leq b} d(n) f(n)$.
Next we deal with $\sum_{a < n \leq b} d(n) \chi(n) f(n)$. Now
\begin{align*}
\sum\limits_{a < n \leq b} d(n) \chi(n) f(n) &= \sum\limits_{a < mn \leq
b} \chi (mn) f(mn)\\
&= \sum\limits_{m \leq b} \sum\limits_{r_{2} = 1}^{k} \chi(mr_{2})
\left( \sum\limits_{\substack{a < mn\leq b\\ n \equiv r_{2} ({\rm mod}\ k)}}
f(mn) \right)\\
&= \sum\limits_{m \leq b} \sum\limits_{r_{2} = 1}^{k} \chi (mr_{2})
\left( \sum\limits_{\substack{a/m < n\leq b/m\\ n \equiv r_{2} ({\rm mod}\ k)}}
f(mn) \right)\\
&= \sum\limits_{m \leq b} \sum\limits_{r_{2} =1}^{k} \chi
(mr_{2}) \left\lbrace \frac{1}{km} \int_{a}^{b} f(u) {\rm d}u\right.\\
&\quad\ \left. + \int_{a}^{b} \psi \left( \frac{(u/m) - r_{2}}{k} \right)
{\rm d}f(u)\right.\\
&\quad\  + f(a) \psi \left( \frac{(a/m) - r_{2}}{k} \right)\\
&\quad\ \left. -f(b) \psi
\left( \frac{(b/m) - r_{2}}{k} \right) \right\rbrace.
\end{align*}
Substituting the series for $\psi$, then interchanging $\sum$ and $\int$
and then integrating by parts in one way or the other as the case may
be, we can get the Euler--Maclaurin or the Poisson analogue of $\sum_{a
< n \leq b} d(n) \chi(n) f(n)$.


\begin{thebibliography}{99}
\bibitem{1} Berndt~B~C, Character analogues of the Poisson and
Euler--Maclaurin summation formulas with application, {\it J. Number
Theory} {\bf 7} (1975) 413--445.

\bibitem{2} Berndt~B~C and Schoenfeld Lowell, Periodic analogues of the
Euler--Maclaurin and Poisson summation formulas with applications to
number theory, {\it Acta Arithmetica} {\bf XXVIII} (1975) 23--67

\bibitem{3} Jutila~M, A method in the theory of exponential sums (Tata
Institute of Fundamental Research, Mumbai) (1987)

\bibitem{4} Voronoi~M~G, Sur une fonction transcendante et sis
applications \`a la sommation de quelques s\'eries, {\it Ann. de l'Ecole
Norm. Sup. (3)} {\bf 21} (1904) 207--267, 459--533
\end{thebibliography}
\end{document}